\documentclass{aims}

\usepackage{txfonts}
\usepackage{amssymb,amsmath,amsthm,mathtools}
\usepackage{graphicx,xparse}
\usepackage{setspace,tikz}
\usepackage{algorithm,algorithmic}

\def\typeofarticle{Research article}
\def\currentvolume{5}
\def\currentissue{x}
\def\currentyear{2020}

\def\ppages{xxx--xxx}
\def\DOI{}
\def\Received{July 2019}
\def\Accepted{September 2019}
\def\Published{December 2019}

\newcommand{\RR}{{\mathbb{R}}}
\newcommand{\ZZ}{{\mathbb{Z}}}

\theoremstyle{definition}

\numberwithin{equation}{section}

\begin{document}

\title{Identifying Transition States of Chemical Kinetic Systems using Network Embedding Techniques}
\author{ Paula Mercurio\affil{1}\corrauth, Di Liu\affil{1}}
  
\shortauthors{the Author(s)}

\address{\addr{\affilnum{1}}
{Department of Mathematics, Michigan State University, East Lansing, MI, USA}
}

\corraddr{Email: mercuri6@msu.edu}

\begin{abstract}
Using random walk sampling methods for feature learning on networks, 
we develop a method for generating low-dimensional node embeddings for directed graphs and 
identifying transition states of stochastic chemical reacting systems. We modified objective functions 
adopted in existing random walk based network embedding methods to handle directed graphs and 
neighbors of different degrees. Through optimization via gradient ascent, we embed the 
weighted graph vertices into a low-dimensional vector space $\RR^d$ while preserving the
neighborhood of each node. We then demonstrate the effectiveness of the method on
dimension reduction through several examples regarding identification of transition states of chemical 
reactions, especially for entropic systems.  
\end{abstract}


\maketitle

\section{Introduction}
Networks provide a natural framework for many chemical and biochemical systems, including chemical 
kinetics and dynamical interactions between biomolecules. The possible stages of a reaction, for example,
can be viewed as vertices in a directed graph. Frequently, real-world networks are high-dimensional and
complex in structure, which leads to difficulties in interpretation and analysis. In order to take full
advantage of the structural information contained in such networks, efficient and robust methods of 
effective network analysis are essential to produce low dimensional representations while retaining 
principal correlations of the graph.

Networks analysis is especially useful for understanding behaviors of nonlinear chemical processes such
as biochemical switches. Biochemical switches are systems of chemical reactions with multiple steady
states. Without random perturbations, the system will spend all its time at the metastable states, where
energy is minimized, while the presence of stochastic noise in the dynamics can lead to transitions
between metastable states. Difficulties arise in the simulation of such systems as a result of the 
multiple time scales at work; transitions between steady states are typically rare, but occur rapidly. 

The different possible states of the reaction, given by the varying concentrations of reactant and product
species, can be viewed as the vertices of a weighted directed graph. The progression of the reaction over
time can then be viewed as a Markov jump process, where the jump probabilities are represented by the
weights on the edges between each pair of nodes. Then, methods for network analysis can be used to
identify interactions and pathways between the components of the biochemical switch. One objective of 
this paper is to suggest a method to analyze this Markov process using node embedding techniques. 
This paper will focus on the $A \to B$ problem: $A$ and $B$ are taken as the reactant and product states,
respectively. The aim of this paper is to investigate the mechanism by which this transition occurs.
    
Recently, methods for analyzing networks through feature representation and graph embedding have
received increasing attention. For an overview on this subject, a number of recent review papers 
are available \cite{survey1, survey2}. While much of this research has been focused on undirected 
graphs, the directed case has also been investigated, for example \cite{cheegerineq,directedreps,zhou}. 
In particular, several methods have been proposed that utilize random walks for node embedding.  
   
Node embedding (or \emph{feature learning}) methods aim to represent each node of a given network 
as a vector in a low-dimensional continuous vector space, in such a way that information about the 
nodes and relationships between them are retained. Specifically, nodes that are similar to each other
according to some measure in the original graph will also have similar representations in the embedding 
space. Usually, this means that the embedding will be chosen in a way that maximizes the inner
products of embedding vectors for embedded nodes that are similar to each other. Compared with the
original complex high-dimensional networks, these low dimensional continuous node representations 
have the benefit of being convenient to work with for downstream machine learning techniques.

Compared with alternative methods that embed edges or entire graphs instead of nodes, node 
embeddings are more adaptive for numerous tasks, including node classification, clustering, and link
prediction \cite{survey2}. Node classification is a process by which class labels are assigned to nodes
based on a small sample of labelled nodes, where similar nodes have the same labels. Clustering
algorithms group the representations of similar nodes together in the target vector space. Link 
prediction seeks to predict missing edges in an incomplete graph.

The recent interest in feature learning has led to a collection of different node embedding methods, 
which can be broadly classified \cite{survey1} by the analytical and computational techniques 
employed such as matrix factorization, random walk sampling, and deep learning network. The 
matrix factorization based methods generate embeddings by first using a matrix to represent the
relationships between nodes, which include the Laplacian matrix as in the Laplacian eigenmaps 
\cite{lapeigs}, the adjacency matrix as in locally linear embedding (LLE) \cite {lle} or graph factorization
(GF) \cite{ahmedetal}, etc. Though the above methods focus on the case of undirected methods,
similar techniques for the directed case have also been suggested \cite{cheegerineq}, where the
development of the Laplacian matrix, a combinatorial Laplacian matrix, and Cheeger inequality allow 
the above approaches to be extended for directed graphs. Factorization based methods such 
as eigenvalue decomposition can apply for certain matrices (e.g. positive semidefinite), 
and additional difficulties may arise with scaling for large data sets. 

Random walk methods are often used for large graphs, or graphs that cannot be observed in their 
entirety. The general idea involves first simulating random walks on a network, then using the output 
to infer relationships between nodes. Random walk methods can be used to study either undirected 
or directed graphs, although much of the previous work has focused on the undirected case. 
For example, DeepWalk \cite{deepwalk} uses short unbiased random walks to find similarities 
between nodes, with node pairs that tend to co-occur in the same random walks having higher similarities. 
This method performs well for detecting \emph{homophily} -- similarity based on node adjacency and 
shared neighbors. Similar methods, such as node2vec \cite{node2vec}, use biased random walks to 
capture similarity based on both homophily and \emph{structural equivalence} -- similarities in the nodes' 
structural roles in the graph.  Both of these methods specifically address undirected graphs. A similar
approach designed for the directed case is proposed in \cite{directedreps}, which uses Markov random 
walks with transition probabilities given by ratios of edge weights to nodes' out degrees, together with the
stationary distribution of the walks, to measure local relationships between nodes. 
    
Theoretical frameworks for the study of transition events of biochemical processes include 
Transition State Theory (TST), Transition Path Sampling (TPS), Transition Path Theory (TPT)
and Forward Flux Sampling (FFS). The main idea of TST is that in order for the system to move from the
reactant state to the product state, the system must pass through a saddle point on the potential energy
surface, known as a Transition State \cite{b:Wigner38}. The TPS technique uses Monte Carlo 
sampling of transition paths to study the full collection of transition paths of a given Markov process
\cite{b:Bolhuis02}. Transition Path Theory (TPT) studies reactive trajectories of a system by considering
statistical properties such as rate functions and probability currents between states \cite{b:EV06}.  
Forward flux sampling \cite{bionetworks}, designed specifically to capture rare switching events in 
biochemical networks, uses a series of interfaces between the initial and final states to calculate rate 
constants and generate transition paths. For the purpose of this paper, we will be using TPT, but 
other such techniques can fit into the method presented here as well.

This paper will suggest a method of using a random walk network embedding approach for node
classification and clustering on directed graphs, as well as identification of transition states 
for the specific case of entropy effects. The method presented reduces the dimension of networks, 
thereby allowing analysis and interpretation of the system. First, in Section 2 we will introduce some
background regarding feature learning on networks, and provide a brief overview of some basic 
principles from TPT. Then Section 3 will introduce a method for node classification for directed graphs,
based on random walk node embedding.  Finally, in Section 4, we will study several examples, showing 
the effectiveness of the method on identifying transition states of entropic systems.

\section{Background}
In this paper, we focus on the mechanism by which a system such as a biochemical switch passes 
through its energy landscape following a \emph{reactive trajectory}. A trajectory is said to be reactive if it
leaves the reactant state $A$ and later arrives at the product state $B$ without first returning to state $A$.
Such a trajectory can be viewed as a sequence of transitions between \emph{metastable states} where 
the energy is locally minimized. In particular, we treat these sequences as Markov jump processes, 
and study the probability space of these reactive trajectories in order to identify and understand the 
transition events of the system. 

Though this paper will be utilizing the techniques of Transition Path Theory for this purpose, it is also 
possible to adopt other frameworks, such as Transition Path sampling, Forward Flux sampling, etc., to
investigate the transition events of a system. As a Monte Carlo technique, Transition path sampling
\cite{b:Bolhuis02} is essentially a generalization of importance sampling in trajectory space. The basic 
idea is to perform a random walk in trajectory space, biased so that important regions are visited more
often.  Forward flux sampling (FFS) \cite{bionetworks}, which generates trajectories from $A$ to $B$ 
using a sequence of nonintersecting parametrizable surfaces, has the advantage that knowledge of 
the phase space density is not required. 

Consider a Markov jump process on a countable state space $S$. 
Let $L=(l_{ij})_{i,j\in S}$ represent the infinitesimal generator of the process. In other words, 
for $i \neq j$, $l_{ij} \Delta t + o(\Delta t)$ represents the probability that a process that is in state $i$ 
at time $t$ will jump to state $j$ during the time interval $[t, t+\Delta t]$. Then the entries of $L$ are
transition rates, satisfying
\begin{align}
\begin{cases}
l_{ij} \geq 0 \hspace{0.1in} &\forall\ i,j \in S,\ i \neq j, \\
\sum_{j \in S} l_{ij} =0, \hspace{0.1in} &\forall\ i \in S.\\
\end{cases}
\end{align}
Let $\{X(t)\}_{t \in \RR}$ represent an equilibrium sample trajectory of the Markov process, 
such that $\{X(t)\}_{t \in \RR}$ is right-continuous with left limits. At time $t$, let the probability 
distribution of this process be denoted by $\mu(t):=(\mathbb{P}(X(t)=i))_{i \in S}^T$. Then the 
time evolution of $\mu(t)$ follows the forward kolmogorov equation (or \emph{master equation})
\[
\frac{d \mu }{dt}=\mu^TL \text{,} \quad t \geq 0.
\]

Denote the time-reversed process by $\{\tilde{X}(t) \}_{t \in \RR}$, and the infinitesimal generator 
of this process by $\tilde{L}$. 
The stationary probability distribution of both processes,   
$\{X(t) \}_{t \in \RR}$ and $\{\tilde{X}(t) \}_{t \in \RR}$ $\pi = (\pi_i)_{i \in S}$, is given by the solution to
\[
0=\pi ^T L\ .
\]

\subsection{Network Embedding for Directed Graphs}

Networks can have large numbers of vertices and complex structures, which can lead to challenges 
in network analysis. To overcome these difficulties and improve understanding of these networks, the
nodes of a given network can be embedded into a low-dimensional vector space, a process called
\emph{feature learning}. The techniques used in this paper to investigate transition states and 
pathways is a combination of neural network-based and random walk-based methods for learning 
latent representations of a network's nodes, as discussed in \cite{node2vec, deepwalk, directedreps}. 
The idea behind these methods is that the nodes in a given network $G=G(V,E)$ can be embedded 
into a low-dimensional vector space in such a way that nodes similar to each other 
according to some measure in the original graph will also have embeddings 
similar to each other in the embedding space. 

Feature learning, or node embedding, requires three main components: an encoder function, 
a definition of node similarity and a loss function. The encoder function is the function that maps 
each node to a vector in the embedding space. To find the optimal embeddings, the encoder 
function can be chosen by minimizing the loss function so that the similarity between nodes 
in the original network corresponds as close as possible to the similarity between their vector
representations in the continuous space. Node similarity can have different interpretations, 
depending on the network and the task to be accomplished. For example,  nodes that are connected, 
share neighbors, or share a common structural role in the original graph can be claimed to be similar, 
and to be embedded in the vector space by representations of close similarity. For the embedded 
vectors, a common approach is to use inner products as the similarity 
measure in the feature space--that is, if $u$ and $v$ are two nodes in a graph with representations 
$z_u$ and $z_v$ in a low-dimensional vector space, then $z_u^Tz_v$ should approximate the similarity
between representations $z_u$ and $z_v$. 

\subsubsection{Factorization Approaches}
Many feature learning methods obtain node embeddings via factorization of a matrix containning 
information regarding relationships between nodes. Locally Linear Embedding (LLE) \cite{lle} and 
Graph Factorization (GF) \cite{ahmedetal}, utilize the node adjacency matrix entry $A_{i,j}$ being the 
weight of edge $(i,j)$ in the graph. Laplacian Eigenmaps and related methods instead factorize the
Laplacian matrix of the whole graph. Optimization of the loss function is typically solve as an eigenvalue
problem.


While the methods mentioned here mainly preserve first-order node relationships, i.e., two nodes that 
are directly connected by an edge are similar, other factorization methods such as GraRep, Cauchy 
graph embedding and structure preserving embedding can detect higher-order relationships between
nodes, or alternative notions of similarity such as structural equivalence \cite{grarep,survey1,survey2}. 

\subsubsection{Random Walk Approaches}\label{sssec:randomwalks}
One technique for finding node embeddings is through the use of random walks. This technique has 
been applied to undirected graphs in methods such as DeepWalk and node2vec 
\cite{deepwalk, node2vec}. In these methods, the general approach is to use short random walks to 
determine similarities for each pair of nodes. The notion of node similarity is co-occurence within a 
random walk; two nodes are similar if there is a higher probability that a walk containing one node will 
also contain the other. The various methods using this approach differ in the implementation of this 
general idea: DeepWalk, for example, uses straightforward unbiased walks, while node2vec uses 
adjustable parameters to bias random walks toward breadth-first or depth-first sampling for a more 
flexible notion of similarity. 

The first stage in the feature learning process is to simulate a certain number of random walks starting 
from each vertex. Using the results of these random walks, the neighbors of each node can be 
determined. The objective is to maximize the probability of observing the neighborhood of each node, 
conditioned on the vector representation of the node. In several of the random walk-based feature 
representation methods for undirected graphs, this is achieved through the Loss function of the form:
\begin{align}
\max_eJ(e)\ =\ \sum_{u \in V} \log Pr(N_S(u)|e(u)),
\end{align}
where $V$ is the set of all vertices in the graph, $N_S(u)$ is the neighborhood of the node $u$, 
and $e(u)$ is the encoder function representating the embedded vector. If we assume that the 
probabilities of observing a particular neighbor are independent, we can have
\begin{align}
Pr(N_S(u)|e(u))=\prod_{n_i \in N_S(u)} Pr(n_i|e(u)).
\end{align}
Since these methods typically aim to maximize the inner products for neighboring nodes, 
they often utilize a softmax function to model each conditional probability:
\[
Pr(n_i|e(u))=\frac{\exp\big(e(n_i) \cdot e(u)\big)}{\sum_{v \in V} \exp\big(e(v) \cdot e(u)\big)}\ .
\]
For networks where $|V|$ is large, techniques such as hierarchical softmax \cite{deepwalk} or 
negative sampling \cite{negsampling} can be used to increase efficiency.
To determine the optimal embedding function $e(u)$, the Loss function can be optimized using gradient
ascent respect to parameters of $e(u)$.

In \cite{directedreps},  a random walk-based method specifically intended for directed graphs is 
proposed. The function to be optimized in this method, when embedding into $\RR^1$,  is 
 \[
I=\sum_u \pi_u \sum_{v, u \to v} p(u,v)(y_u-y_v)^2,
\]
where $y_u$ is the one-dimensional embedding of $u$,
$\pi_u$ is the stationary probability for the node $u$ and $p(u,v)$ represents the random walk's 
transition probability from $u$ to $v$. The transition probabilities are calculated according to
\[
p(u,v)=w(u,v)/\sum_{k, u \to k} w(u,k),
\]
where $w(u,v)$ is the weight on the edge $(u,v)$.
Note that this method, though it uses transition probabilities of a random walk, does not actually 
require simulation of random walks, since all the summations can be done deterministically.

According to the spectral graph theory \cite{cheegerineq}, a circulation on a directed graph is 
a function $F$ such that, for each $u \in V$, $\sum_{u, u \to v} F(u,v)=\sum_{w, v \to w} F(v,w)$.
The Laplacian of a directed graph is therefore defined as 
\[
\mathcal{L}=I-\frac{\Phi^{1/2}P \Phi^{-1/2} + \Phi^{-1/2}P \Phi^{1/2}}{2},
\]
where $\Phi$ is a diagonal matrix with diagonal entries given by 
$\Phi(v,v)=\phi(v)=\sum_{u, u \to v} F(u,v)$ for some circulation function $F$.  
The following definition of the \emph{combinatorial Laplacian} is also due to Chung\cite{cheegerineq}:
\[
\mathbb{L}=\Phi-\frac{\Phi P +P^T \Phi}{2},
\]
where $P$ is the transition matrix, i.e. $P_{ij}=p(i,j)$, $\Phi$ is the diagonal matrix of the stationary
distribution, i.e. $\Phi=diag(\pi_1,\dots,\pi_n)$. Note that the combinatorial Laplacian is symmetric and
semi-positive definite. It is shown in \cite{directedreps} that
\[
\sum_{u} \pi_u \sum_{v, u \to v} p(u,v)(y_u-y_v)^2=2y^T\mathbb{L}y
\]
where $y=(y_1, y_2, \dots , y_n)^T$. This result demonstrates that this method is analogous to the
Laplacian-based methods for undirected graphs, e.g. Laplacian eigenmaps \cite{survey1}.

Random walk methods have been shown to perform well compared to other methods, and can be 
useful for a variety of different notions of node similarity. Previous work has shown them to be robust, 
efficient, and capable of completing a diverse range of tasks, including node classification, 
link prediction, clustering, and more \cite{survey1, survey2, deepwalk, node2vec}. Because they do 
not examine the entire graph at once, these methods are also useful for very large networks and 
networks that cannot be observed in their entirety. 

\subsubsection{Neural Network Approaches}\label{sssec:neuralnets}
A third class of graph embedding methods involves the use of deep learning techniques, including 
neural networks \cite{survey2,survey1}, e.g. Graph Neural Network (GNN) \cite{survey2, gnn}.
Such methods have been shown to perform well, particularly for dimension reduction and identifying
nonlinear relationships among data. Random walk based methods, such as DeepWalk, can incorporate
deep learning algorithms like SkipGram for graph embedding, where random walks serve to provide 
neighborhood information as input. 

In general, neural network methods typically work by assigning a weight to each node representation, 
and then combining those terms via a transfer function. This result is then used as the input for an
activation function, such as a sigmoid function, i.e. $\sigma(x)=1/(1+\exp(-x))$. This process may 
then repeat for multiple layers. In a neural network method, the parameters to be optimized are the 
weights, which are updated after each layer, typically by gradient descent. 
In the GNN \cite{gnn}, the objective function to be minimized is of the form
\[
W=\sum_i \sum_j (t_{i,j}-\varphi_w(G_i,n_{i,j}))^2,
\]
where $G_i$ is the learning set, the $n_{i,j}$ and $t_{i,j}$ are the nodes and target outputs, and the 
goal is to choose the parameter $w$ in such a way that $\varphi$ closely approximates the target 
outputs. This optimization occurs through gradient descent, where the gradient is computed with 
respect to the weights $w$. 

\subsection{Transition Path Theory}\label{subsec:tpt}
In order to identify the transition paths of a given Markov process, we will adopt the framework of
Transition Path Theory (TPT). The following notations are mostly from \cite{tpt}. Consider a state 
space $S$, for an initial state $A \in S$ and final state $B \in S$, a trajectory $X(t)$ is said to be 
{reactive} if it begins from state $A$ and reaches at state $B$ before returning to state $A$.

To determine whether a given reaction path is reactive, we will need the following forward 
and backward committor functions.
For each $i \in S$, the forward committor $q_i^+$ is the probability that a process initially 
at state $i$ will reach state $B$ before it reaches state $A$. Similarly, the discrete backward committor 
$q_i^-$ denotes the probability that a process arriving at state $i$ was more recently in state $A$  
rather than in state $B$.  The forward committors satisfy the following Dirichlet equation
\begin{align*}
\begin{cases}
\sum_{j \in S} l_{ij}q_j^+ &=0 \hspace{0.5in} \forall i\in (A \cup B)^C \\
q_i^+&=0 \hspace{0.5in} \forall i \in A \\
q_i^+&=1 \hspace{0.5in} \forall i \in B,\\
\end{cases}
\end{align*} 
and similar equations are satisfied by backward committors.

The probability current, or flux, of reactive trajectories gives the average rate at which reactive paths 
transition from state to state. For trajectories from $A$ to $B$, the probability current is defined for 
all $i \neq j$ such that
\begin{align*}
f_{ij}^{AB}&=\lim_{s \to 0+} \frac{1}{s} \lim_{T \to \infty} 
\frac{1}{2T} \int_{-T}^T \mathbf{1}_{\{i\}} (X(t))\mathbf{1}_{\{j\}} (X(t+s)) 
\times \sum_{n \in \ZZ} \mathbf{1}_{(-\infty , t_n^B]} (t) \mathbf{1}_{[t_n^A, \infty)} (t+s) dt \\ &=f^{AB}_{ij}.
\end{align*}
Also, for all $i \in S$, $f^{AB}_{ii}=0$. It can be shown that
\begin{equation*}
f^{AB}_{ij}=
\left\{\begin{array}{lll}
\pi_iq_i^-l_{ij}q_j^+, \quad &\text{if } i \neq j \\
0  &\text{if } i = j. \\
\end{array}\right.
\end{equation*}

Since we are primarily interested in the flow from $A$ to $B$, and the process can move in either 
direction between two adjacent nodes on the path, the following {effective current} accounts for the 
net average number of reactive paths that jump from $i$ to $j$ per time unit:
\[
f_{ij}^+ = \max(f_{ij}^{AB}-f_{ji}^{AB}, 0).
\]
We can also define the total effective current for a node $i$:
\[
C_i^+=\sum_{j:(i,j)\in E} f_{ij}^+.
\]
In \cite{du}, the effective current is generalized to a connected subnetwork of nodes $\Omega$ such that
\[
C^+(\Omega)=\sum_{i \in \partial \Omega} C_i^+,
\]
which leads to the following definition of transition states of a time reversible process as subsets of the 
state space:
\[
TS=\lim_{\sigma \to 0} \underset{\Omega}{\mathrm{argmax}} 
\{ C^+(\Omega) \exp\left(-(q^+-0.5)^2/\sigma^2\right) \}
\]
For the non-time reversible case, transition states can be identified similarly, replacing the exponential
argument $\left(-(q^+-0.5)^2/\sigma^2\right)$ with $(-(q^+-0.5)^2+(q^--0.5)^2)/\sigma^2$.
The above definition gives the flexibility of identifying transition states as subsets of the states space, 
which will be very useful when dealing with complex dynamics such as entropy effects. 

\section{Identifying Transition States using Network Embedding}
We now introduce a method to identify transition states and paths of Markov processes, by random walk 
based Network Embedding techniques for directed networks. Consider a Markov jump process in a 
countable state space $S$. Using the infinitesimal generator $L$ for the process, we can calculate the 
stationary probability distribution $\pi$ and the forward and backward committors, based on which 
probability current and effective currents $f_{ij}^{AB}$ and $f_{ij}^+$ can then be computed.
Once these quantities have been obtained, the discretized state space can be viewed as a weighted, 
directed graph $G(V,E)$. The nodes $v \in V$ are the grid points representing states of the system, 
and each pair of adjacent nodes are connected by an edge $(u,v) \in E$ with weights given by the 
effective probability current, such that the direction of the edge will be determined by the sign of the 
effective current. 

Once $G$ is constructed, we can apply feature learning techniques for node classification to identify 
transition states. 
Random walk trials for similarities between nodes will start from each node in the space, using the edge
weights as transition probabilities. The outputs of the random walks can then be used to compute
experimental ``neighbor probabilities", i.e., the probability that a random walk of length $d$ starting at 
node $m$ will contain node $n$. If the neighbor probability for a node pair $(n,m)$  is sufficiently high, 
then the node $n$ is claimed as a neighbor of $m$. Note that this is a directed process; $n$ can be a 
neighbor of $m$ even if $m$ is not a neighbor of $n$.

After finding these conditional probabilities and identifying the neighborhoods of each node, the next 
step is to maximize an objective function in order to solve for the node embeddings.
As in the node2vec and DeepWalk methods for the undirected case, the conditional probabilities are 
modeled using softmax functions. In order to adapt this technique for the case of directed graphs, we
will use a modified version of the objective function: 
\begin{equation}\label{eq:objnd}
V[e]=\sum_{u} \pi_u \sum_{v \in N(u)} Pr(v|e(u)) .
\end{equation}
Optimizing this function will give the encoder function $e^{opt}$ that maximizes the probability of
the neighborhood of each node, with the most relevant nodes having the most influence on the feature 
representations.  

In the original form of the objective function \eqref{eq:objnd} for the undirected graph methods described 
above, the softmax functions depend only on whether nodes are neighbors, and is independent of the 
probability (frequency) that a particular neighbor will show up in a random walk trial. In other words, all 
neighbors are treated equally, while in reality some node pairs are more likely to co-occur than others. To 
address this, the objective function used here incorporates the neighboring probabilities, denoted
$NP(\cdot, \cdot)$ such that 
\[
Pr(n_i|e(u))=\frac{\exp(e(n_i) \cdot e(u) )NP(n_i,u)}{\sum_{v \in V} \exp(e(v) \cdot e(u))NP(v,u)} .
\]
\begin{algorithm}[h]
\caption{A random walk simulation method for directed graphs, given the transition matrix $t$.}
\begin{algorithmic}
\FOR{$u \in V$}
\STATE $walk(1)=u$
\FOR{ $step=1:walklength$}
\FOR{ $k=2:d$}
\STATE $total=\sum_{v, walk(k-1) \to v} t(walk(k-1),v)$
\FOR{ $v : walk(k-1) \to v$} 
\STATE $P(walk(k-1),v)=t(walk(k-1),v)/total$
\ENDFOR
\STATE Use probability $P(walk(k-1), \cdot)$ to choose $walk(k)=v$ for some neighbor of 
$walk(k-1)$.\\ 
Set $counter(walk(k),u)=counter(walk(k),u)+1$.
\ENDFOR
\ENDFOR
\ENDFOR
\FOR{$u,v \in V$}
\STATE $neighborprobability(u,v)=counter(u,v)/\sum_k counter(u,k)$
\ENDFOR
\end{algorithmic}
\end{algorithm}

%

Assuming $e(u)=E_2u$, for $E_2$ being a $2 \times d$ matrix, the objective defined above can then be
maximized using gradient ascent to obtain the matrix $E_{2}^{opt}$ and therefore the feature embeddings. 
In practice, the linear assumption on the map from nodes to embeddings appears to be insufficiently
flexible, particularly for higher-dimensional problems. To improve upon this, a neural network can be 
incorporated into the feature learning process, with the sigmoid function $\sigma(u)=1/(1+\exp(-e(u)))$ 
where
$e(u)$ is the representation in the embedding space of the node $u$. The current algorithm uses a
relatively small number of layers; raising this number results in node representations that are clustered
more closely to their neighbors. 

To obtain transition states for a particular process, we will examine the similarities between the node 
representations of each state and representation of metastable state $A$. Here, by ``similarities" we mean
the conditional probabilities given by the softmax units discussed previously.
By construction, these values will naturally be small for nodes that are further from node $A$, due to 
hopping distances that are longer than the length of the random walks. So, to find the transition states, we 
combine similarities of node pairs with shorter hopping distances via matrix propagation, e.g. 
\begin{equation}\label{e:prop}
sim(A,u)=\sum_{ A \to v, v \to u} sim(A,v)*sim(v,u).
\end{equation}
Note that a node $u$ for which this value is high is
a higher-order neighbor of the reactant node $A$; therefore, such a node will have a higher probability of 
lying along a transition path.
We then define transition states as nodes with high probability of lying on a transition path, which are not 
direct neighbors of the reactant or product states.

\section{Examples and Results}
Next, we will illustrate the suggested method on several simple examples.
In the examples presented here, we have used $100$ random walk trials, each of length $9$ steps. 
Smaller numbers of trials increase the impact of noise in the results. Due to the matrix propagation
technique described above, changes to the walk length parameter do not have a noticeable affect.

\subsection{Diffusion Process with Energy Barrier}

\begin{figure}[t] 
\begin{center}
\includegraphics[scale = 0.4]{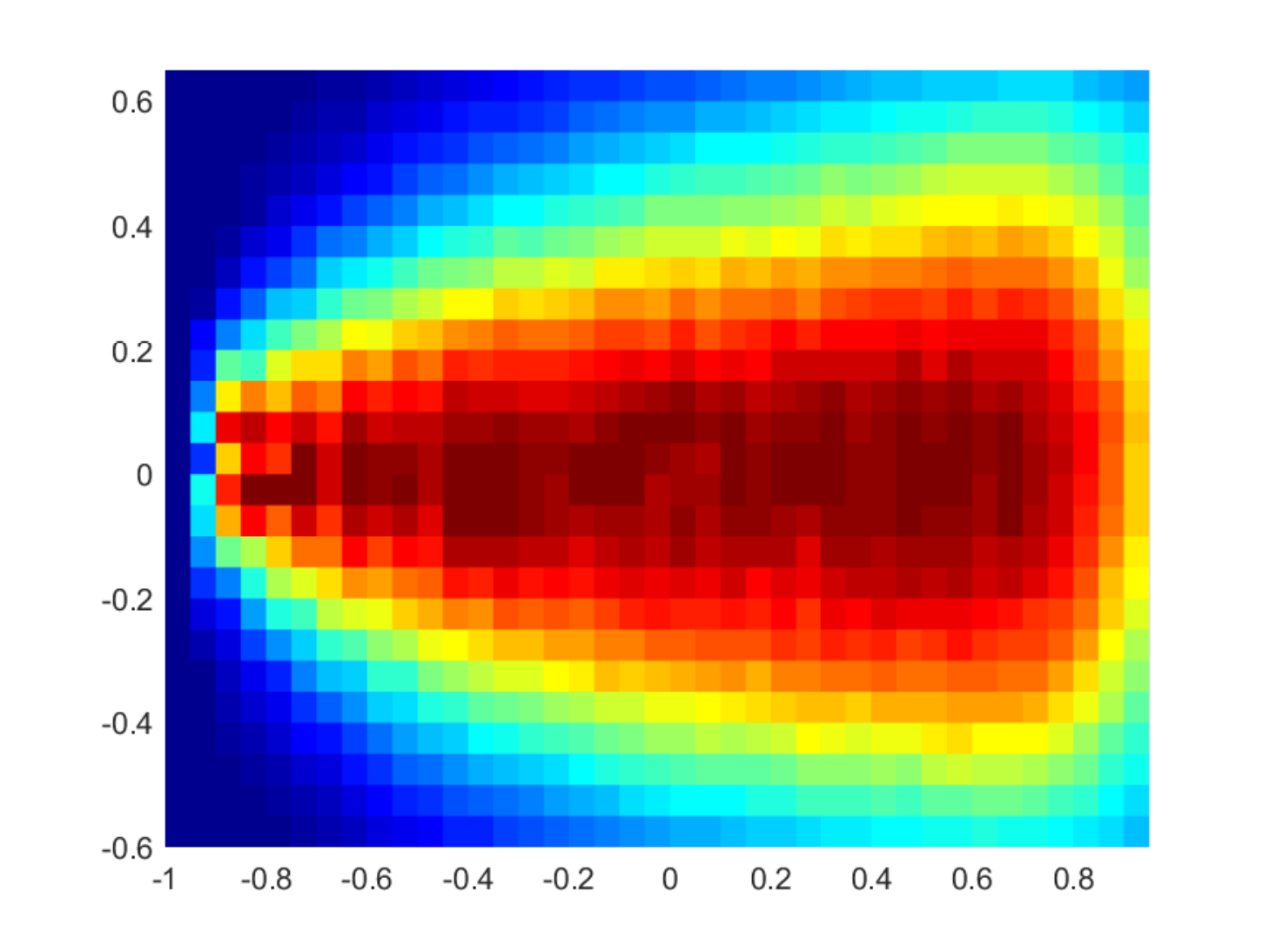}%
\label{fig1}%
\hfil
\includegraphics[scale = 0.4]{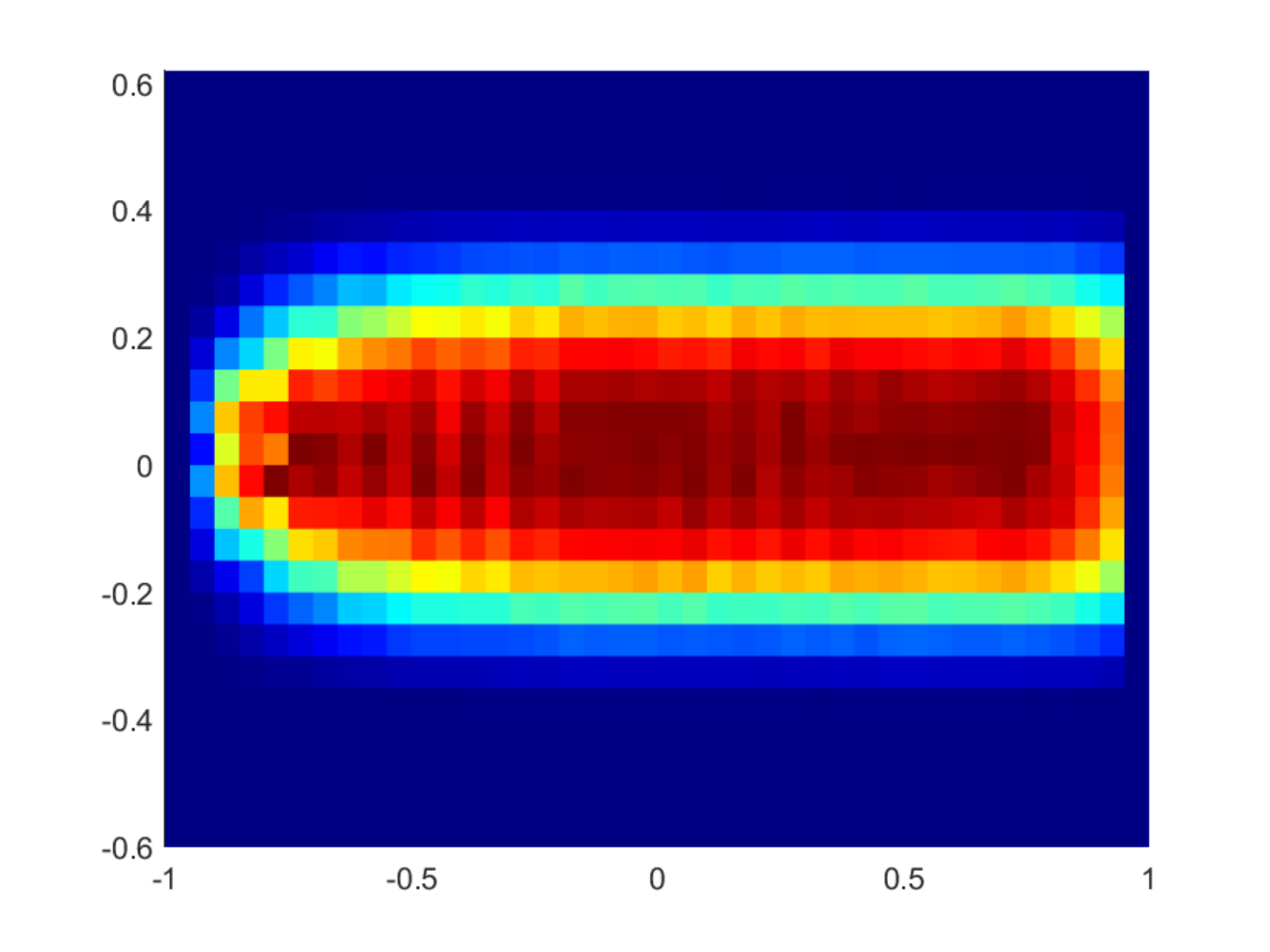}
\caption{Diffusion process with energy barrier:
plot assigning colors to each node in the discretized domain $\Omega$, according to each node's 
similarity to the starting node $A$. 
Higher similarity values correspond to higher dot products between the 
two-dimensional embeddings of the nodes. The left plot corresponds to a parameter value of 
$\varepsilon=0.01$, while the right plot corresponds to $\varepsilon=1$.\label{fig:diff}}%
\end{center}
\end{figure}

\begin{figure}[t] 
\begin{center}
\includegraphics[scale = 0.6]{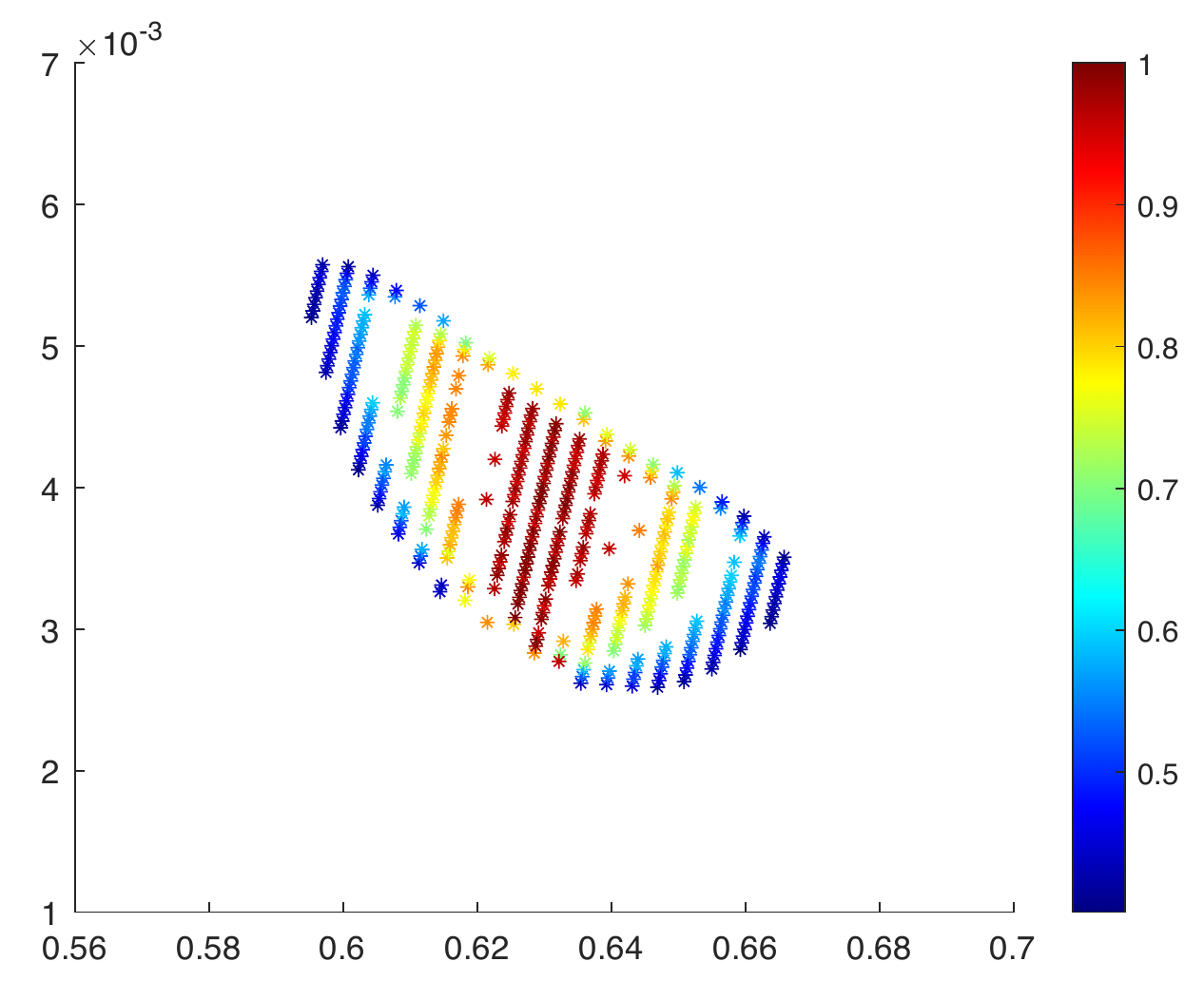}%
\caption{Feature representations of nodes with high similarity to the reactant state $A$ for the 2D diffusion
process with energy barrier, with $\varepsilon=0.01$. \label{fig:2ddiffgroups}}%
\end{center}
\end{figure}

First we will consider a two-dimensional diffusion process representing the position of a particle which
follows the SDE  
\begin{align*}
\frac{dx}{dt}&=-\frac{dV}{dx}+\frac{dW_1(t)}{dt}, \\
\frac{dy}{dt}&=-\frac{dV}{dy}+\frac{dW_2(t)}{dt},\\
\end{align*}
where $V$ is the potential
\begin{align}\label{eq:V}
    V = (x^2-1)^2 + \varepsilon y^4,
\end{align}
and $W_1(t)$ and $W_2(t)$ are independent Brownian motion processes.

This potential $V$ has two metastable states, or local minima, at (-1,0) and (1,0), as well as a saddle 
point at the origin. Here I will take (-1,0) to be the reactant state $A$, and (1,0) to be the product state $B$.
The diffusion process can be approximated using a Markov (birth-death) jump process on a discrete state
space. For this example, we will study the diffusion process on the domain $\Omega = [-1,1]\times
[-0.75,0.75]$ by examining  the Markov jump process on a grid 
$D=((-1+h\ZZ) \times (-0.75+h\ZZ)) \cap \Omega$, with $h=0.05$.

In order to examine the behavior of this Markov process, we should first compute the probability currents 
for the process. The infinitesimal transition matrix $L$ can be found by using the jump rates for the 
birth-death process $l_{nm}$ for each pair of adjacent nodes $n,m$ \cite{handbook}. 
We define the following constants for a state $(x,y) \in (a,b)$:
\begin{align*}
k_x^+(x,y) &=  1/(2*h^2)-1/(2*h)*dVdx(x-h), \\
k_x^-(x,y)&= 1/(2*h^2)+1/(2*h)*dVdx(x+h),\\
k_y^+ (x,y)&= 1/(2*h^2)-1/(2*h)*dVdy(y-h), \\
k_y^- (x,y)&= 1/(2*h^2)+1/(2*h)*dVdy(y+h). 
\end{align*} 
Also let $k_x^+(x+h,y)=1/h$ if $x= A$, $\ k_x^+(x+h,y)=0$ if $x =B$, and 
$k_x^-(x-h,y)=0,\ k_x^-(x-h,y)=1/h$ if $x =A$ or $x =B$ respectively, and similarly for $k_y^+$ and
$k_y^-$.
Then the infinitesimal generator can be expressed in terms of its action on a test function $f$ such that
\begin{equation*}
\begin{aligned}
    (Lf)(x,y)=&k_x^+(x+h,y)(f(x+h,y)-f(x,y)) 
    + k_x^-(x-h,y)(f(x-h,y)-f(x,y)) \\
    &+ k_y^+(x,y+h)(f(x,y+h)-f(x,y)) 
    + k_y^-(x,y-h)(f(x,y-h)-f(x,y)). 
\end{aligned}
\end{equation*}

Figure \ref{fig:diff} assigns each point in the grid a color based on its similarity to the node of  $A=(-1,0)$,
for epsilon values $\varepsilon=.01$ and $\varepsilon=1$ in equation \eqref{eq:V}, respectively. To 
account for the fact that nodes with greater hopping distances from $A$ will have lower similarities to this 
node, we propagate the similarity matrix according to \eqref{e:prop} in order to assign similarities to more
distant nodes. The red area passing from $A$ to $B$ represents the region that reactive trajectories will 
pass through with the highest probability for each value of epsilon.

In Figure \ref{fig:2ddiffgroups}, the node embeddings $(u,v)=e_{opt}(x,y)$ corresponding to a subset of nodes 
with higher similarity values to $A$ are shown, for the case where $\varepsilon=0.01$, where $e_{opt}$ is the 
optimal linear encoding. Notice that in this figure, the embeddings for nodes with the highest probability of 
occurring in a reactive trajectory are clustered together, while nodes with lower probabilities tend to be grouped 
with other nodes with similar probabilities; more specifically, nodes with similarity values in the range $0.5-0.6$ 
have embeddings clustered near the upper left or lower right of the figure, while nodes with values from $0.7-0.8$
will have embeddings which appear in one of the two orange clusters.

In the context of the diffusion process with the potential $V$, the \emph{entropy effect} refers to changes in 
the system's observed dynamics in response to the change of the parameter $\varepsilon$. That is, as 
$\varepsilon$ is decreased, $dV/dy$  shrinks and therefore the negative term of $dy/dt$ becomes small. 
As a result, the shape and size of the saddle point located at the origin is altered, and we can expect the 
y-coordinate to take on a greater variety of values with higher probabilities in this case.
In particular, around the origin, smaller values of $\varepsilon$ should result in higher similarities between the 
nodes above and below the origin, since in this case the process is more likely to move up and down. For 
larger values of $\varepsilon$, the process is expected to move straight forward from $A$ to $B$, with a 
smaller probability of moving up or down. 

\begin{figure}[t] 
\begin{center}
\includegraphics[scale = 0.6]{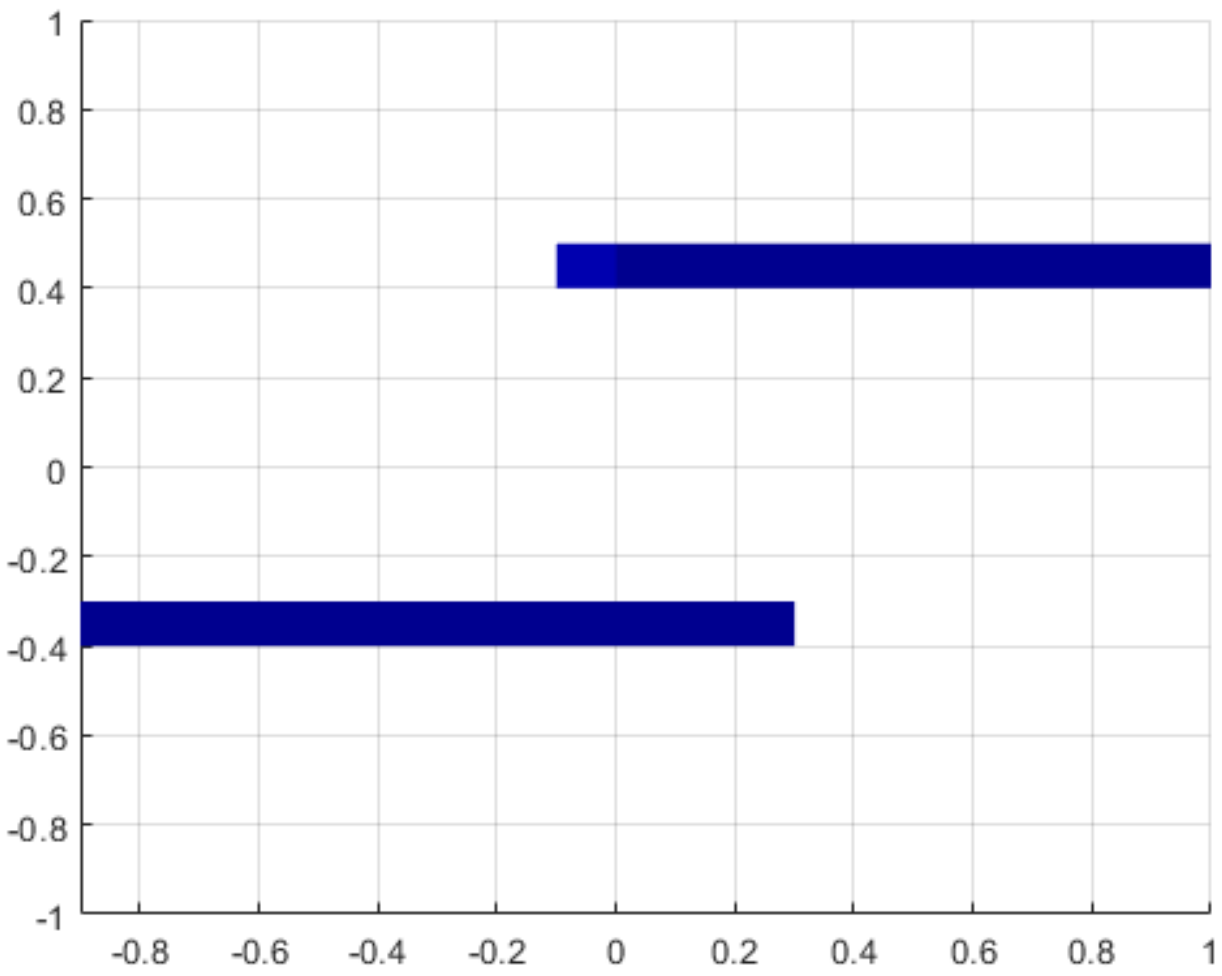}
\qquad
\includegraphics[scale = 0.6]{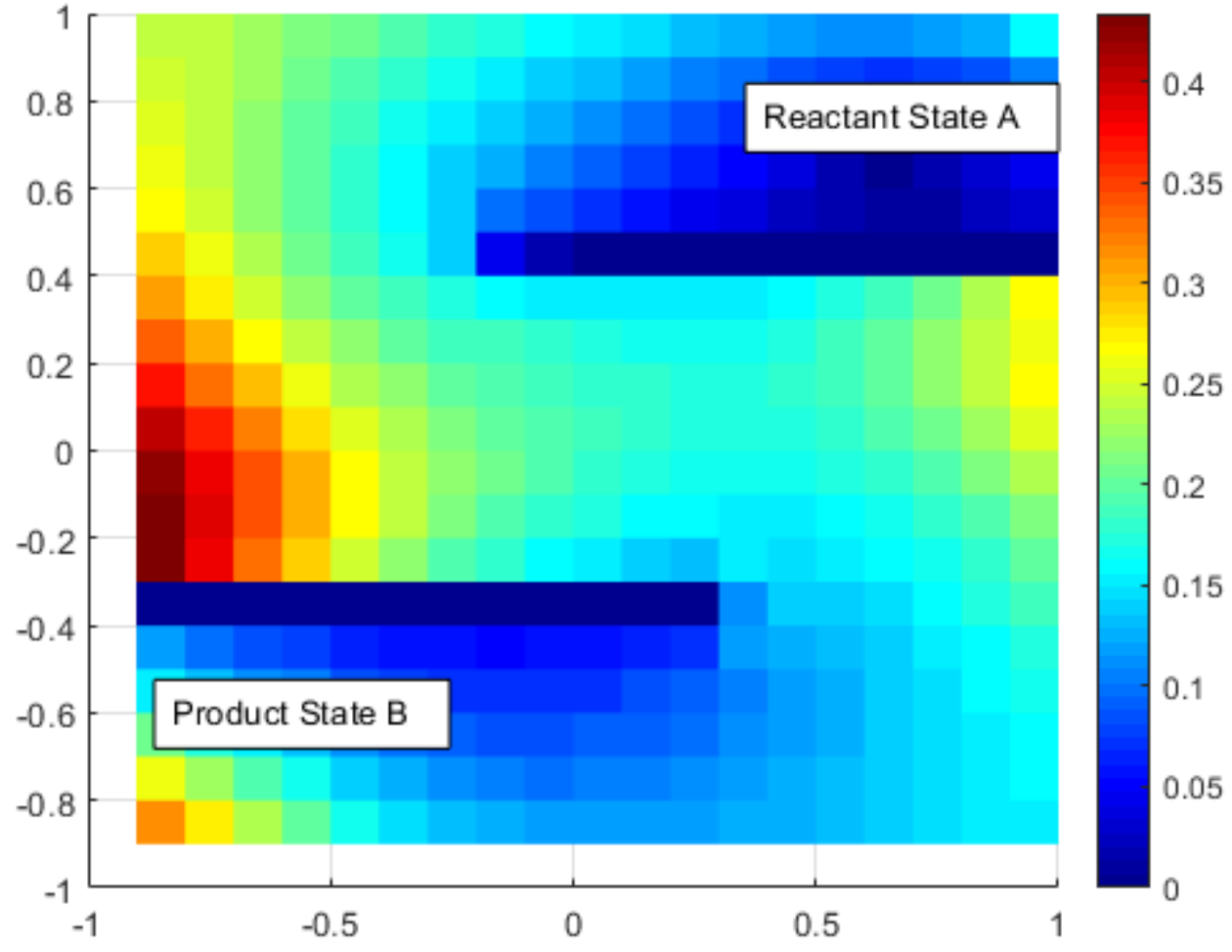}%
\caption{Diffusion process with entropic bottleneck: surface plot of the similarity function between the reactant 
state $A$ and nodes in $\Omega$. Here, the discretization with $h=0.1$ is used. \label{fig:walls}}%
\end{center}
\end{figure}
\subsection{Diffusion Process with Barriers}
\begin{figure}[h] 
\begin{center}
\includegraphics[scale = 0.6]{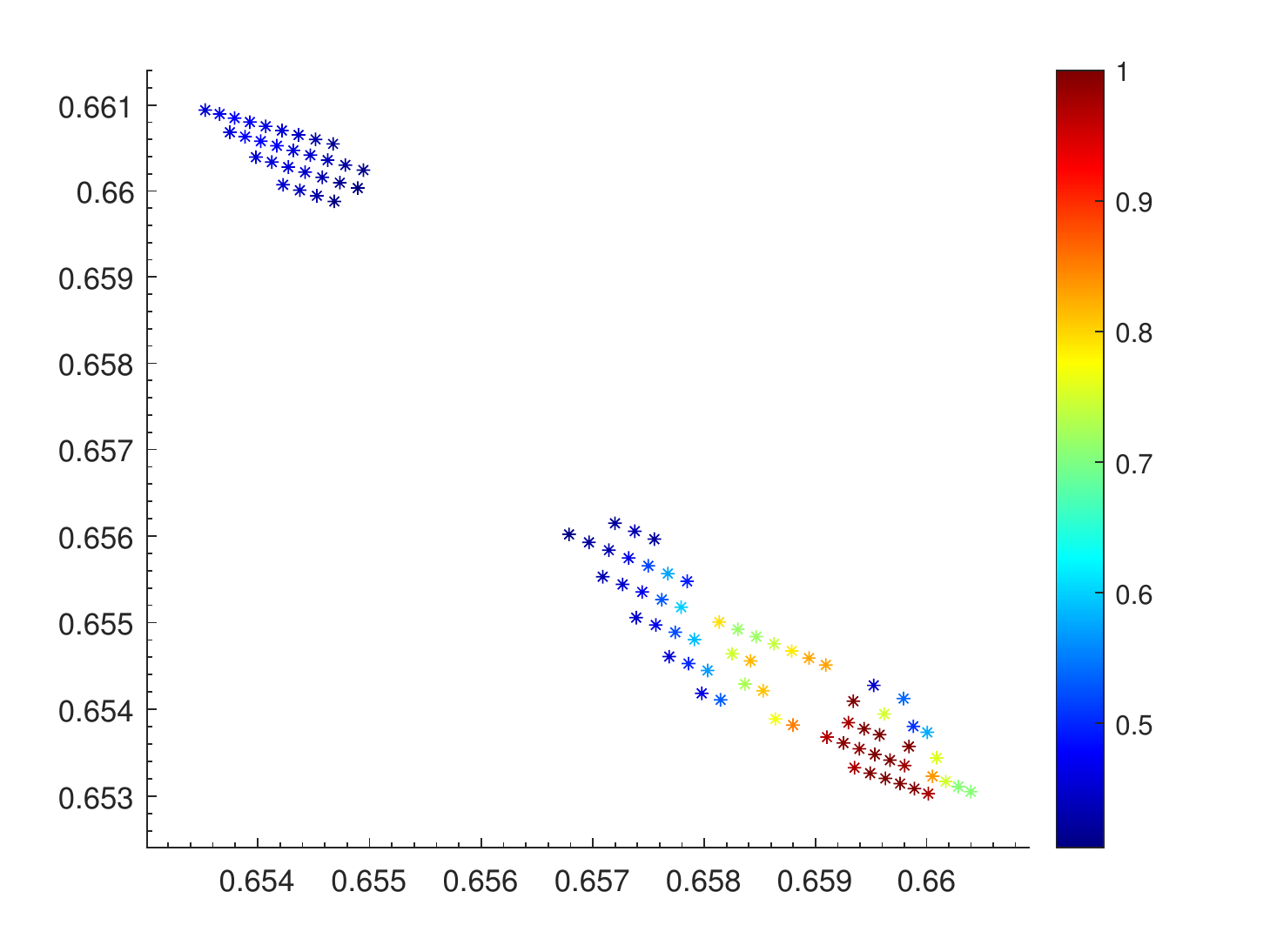}%
\caption{Feature representations of nodes with high similarity to the reactant state $A$ for the entropic 
diffusion process. \label{fig:purediffgroups}}%
\end{center}
\end{figure}

Next we consider a pure diffusion process ($V(x,y)=0$) on the domain $\Omega=[-1,1]\times [-1,1]$ with two 
barriers as depicted in Figure 3. This example was also discussed in \cite{tptexamples}. Here we take the node at 
$(0.6, 0.6)$ to be the initial state $A$, and $(-0.6, -0.6)$ as the final state $B$. The particle spends most of its time
between the reactant and product states, in the region between the barriers ($[-1,1]\times [-0.4,0.4]$). 
The only factor affecting the probability that the particle will reach $B$ before returning to $A$ is its current 
distance to $B$, since there are no energetic obstacles to overcome. In order to travel from $A$ to $B$, the 
process must find its way past the obstacles by chance, thereby overcoming an \emph{entropic barrier}.  

Figure 4 shows the node embeddings for a subset of nodes with high similarity values to reactant $A$. Note that 
the color scheme used in this figure corresponds to the color scheme in Figure 3; node embeddings that are 
colored red in the scatter plot correspond to red areas in the surface plot. In Figure 4, node
embeddings again tend to be grouped by the corresponding values of the similarity function: nodes that are very
similar to the initial state $A$ (similarity function value $>0.9$) appear in the red cluster toward the lower right. 
Other clusters of different colors correspond to groups of nodes with lower probabilities of appearing in a reactive 
trajectory. The yellow cluster contains embeddings of nodes with similarity function values approximately $0.7-0.8$,
while the blue cluster corresponds to nodes with similarity values near $0.5-0.6$. 

\subsection{Three-Dimensional Toggle Switch}
Now we will consider a higher-dimension example, a stochastic model of a 3D genetic toggle switch consisting 
of three genes that each inhibit the others' expression \cite{du}. We consider the production and 
degradation of the three gene products, $S_1, S_2,$ and $S_3$:
\[
* \underset{\alpha_4}{\stackrel{\alpha_1}{\rightleftharpoons}} S_1, * \underset{\alpha_5}{\stackrel{\alpha_2}
{\rightleftharpoons}}S_2, * \underset{\alpha_6}{\stackrel{\alpha_3}{\rightleftharpoons}} S_3,
\]
where the parameters $\alpha_i$ are defined:
\begin{align*}
&\alpha_1 = \frac{c_{11}}{(65+x_2^2)(65+x_3^2)}\ , 
&&\alpha_2 =  \frac{c_{12}}{(65+x_1^2)(65+x_3^2)}\ , \\
&\alpha_3 =  \frac{c_{13}}{(65+x_1^2)(65+x_2^2)}\ , 
&&\alpha_4 = c_4x_1\ ,\\
&\alpha_5 = c_5x_2, 
&&\alpha_6 = c_6x_3, 
\end{align*}
with $c_{11}=2112.5$, $c_{12}=845$, $c_{13}=4225$, $c_4=0.0125$, $c_5=0.005$, and $c_6=0.025$.
From the stationary probability distribution, it can be seen that this model has three metastable states:
\begin{align*}
A=\{x \in S &| 35 \leq x_1 \leq 45, 0\leq x_2 \leq 4, 0 \leq x_3 \leq 4 \}, \\
B=\{x \in S &| 35 \leq x_2 \leq 45, 0\leq x_1 \leq 4, 0 \leq x_3 \leq 4 \}, \\
C=\{x \in S &| 35 \leq x_3 \leq 45, 0\leq x_2 \leq 4, 0 \leq x_1 \leq 4 \}. 
\end{align*}
Here we will choose $A$ to be the reactant state and $B$ to be the product state. 

\begin{figure}[t] 
\begin{center}
\includegraphics[scale = 0.6]{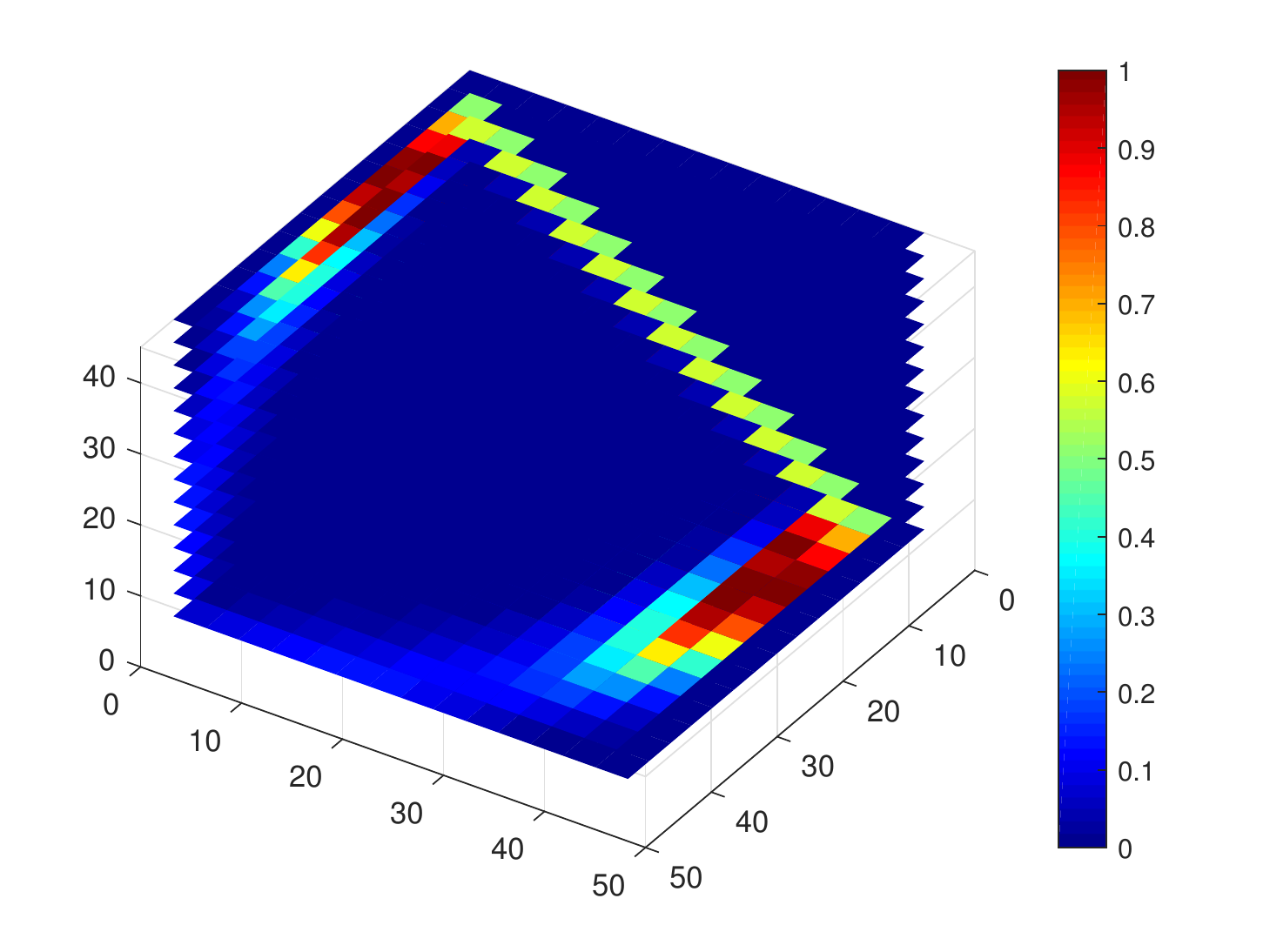}%
\caption{Surface plot of the similarities to $A$ of grid points in the domain $\Omega$ for the 3-D
Toggle switch model. Results are for a $15 \times 15 \times 15$ discretization ($h=3$). \label{fig:3dtoggle}}%
\end{center}
\end{figure}

Figure 5 assigns colors to each node in 
$\Omega = \{ih,jh,kh|i,j,k \in \ZZ \} \cap [0,45] \times [0, 45] \times [0,45]$ based on their similarity to node sin $A$,
propagating the similarity matrix for nodes with larger hopping distance from $A$ as described previously.
This figure shows two potential transition states for this system. A reactive trajectory from $A$ to $B$ will be most
likely to pass directly from $A$ to $B$ along the higher-probability trajectory, represented in red. However, such a 
trajectory may instead reach $B$ after travelling through state $C$, as indicated by the fainter region connecting 
$A$ and $B$ via $C$. 

Below, Figure 6 shows the two-dimensional node embeddings in $\RR^2$, with the embeddings of nodes with zero
similarity to $A$ omitted for clarity. There are several distinct clusters present in Figure 6, each representing a 
cluster of similar nodes. The nodes with the highest similarities to $A$ after propagation of the similarity matrix 
(represented in red) are shown in the cluster to the lower left of the figure. Hence, the nodes in this cluster represent
the nodes expected to appear in the transition path. 

\begin{figure}[h] 
\begin{center}
\includegraphics[scale = 0.6]{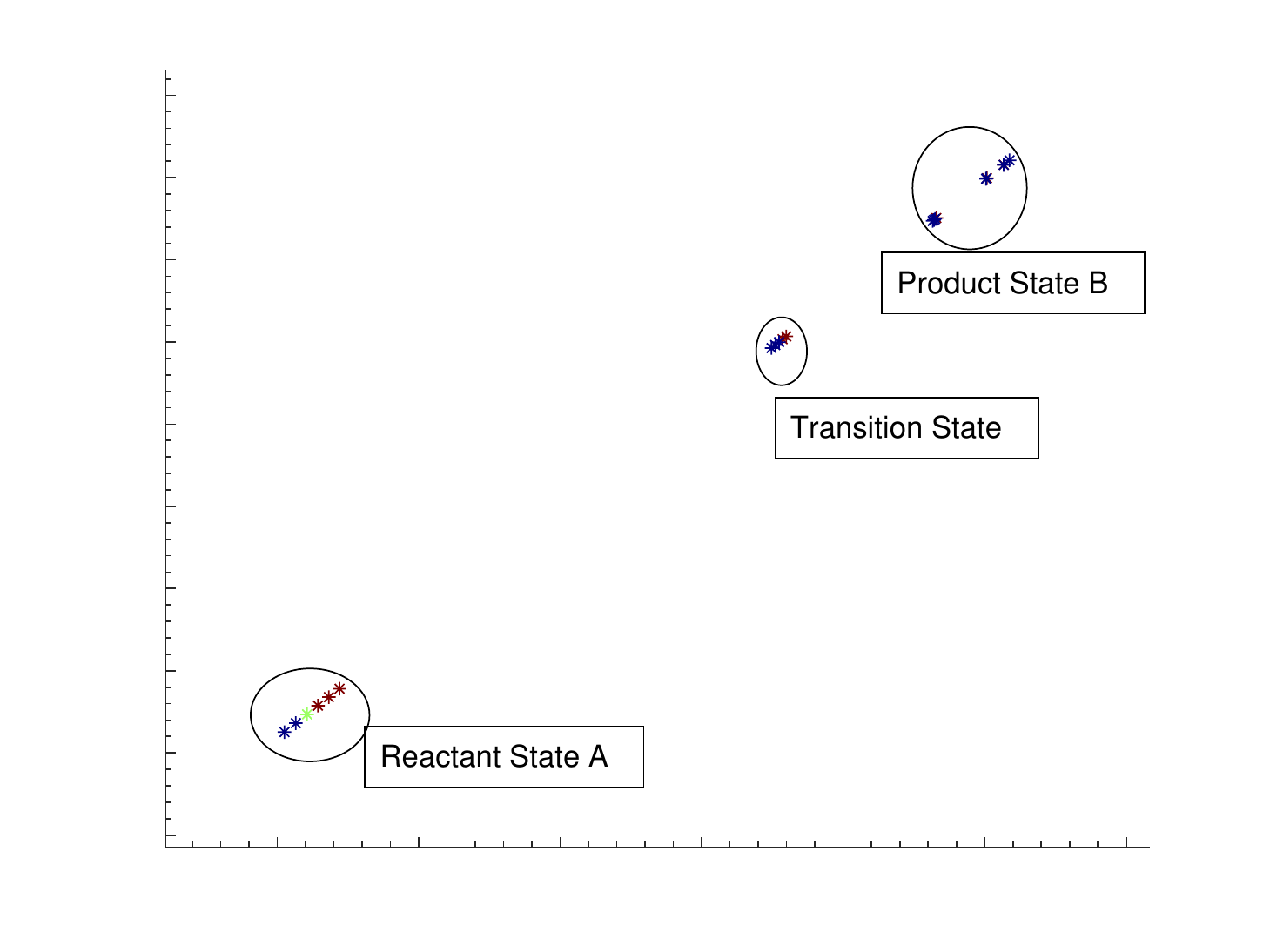}%
\caption{Feature representations of nodes with high similarity to the reactant state $A$ for the 3D toggle switch.
\label{fig:3dtogglegroups}}%
\end{center}
\end{figure}

\subsection{Stochastic Virus Model}
Now we will consider a higher-dimension example, a 3D stochastic model for virus propagation consisting of 
three species which are necessary for virus production \cite{stochvirus}.  The three species involved in this 
system are the template (\emph{tem}),  viral genome (\emph{gen}), and structural (\emph{struct}) proteins of 
the virus. 

\begin{figure}[t] 
\begin{center}
\includegraphics[scale = 0.6]{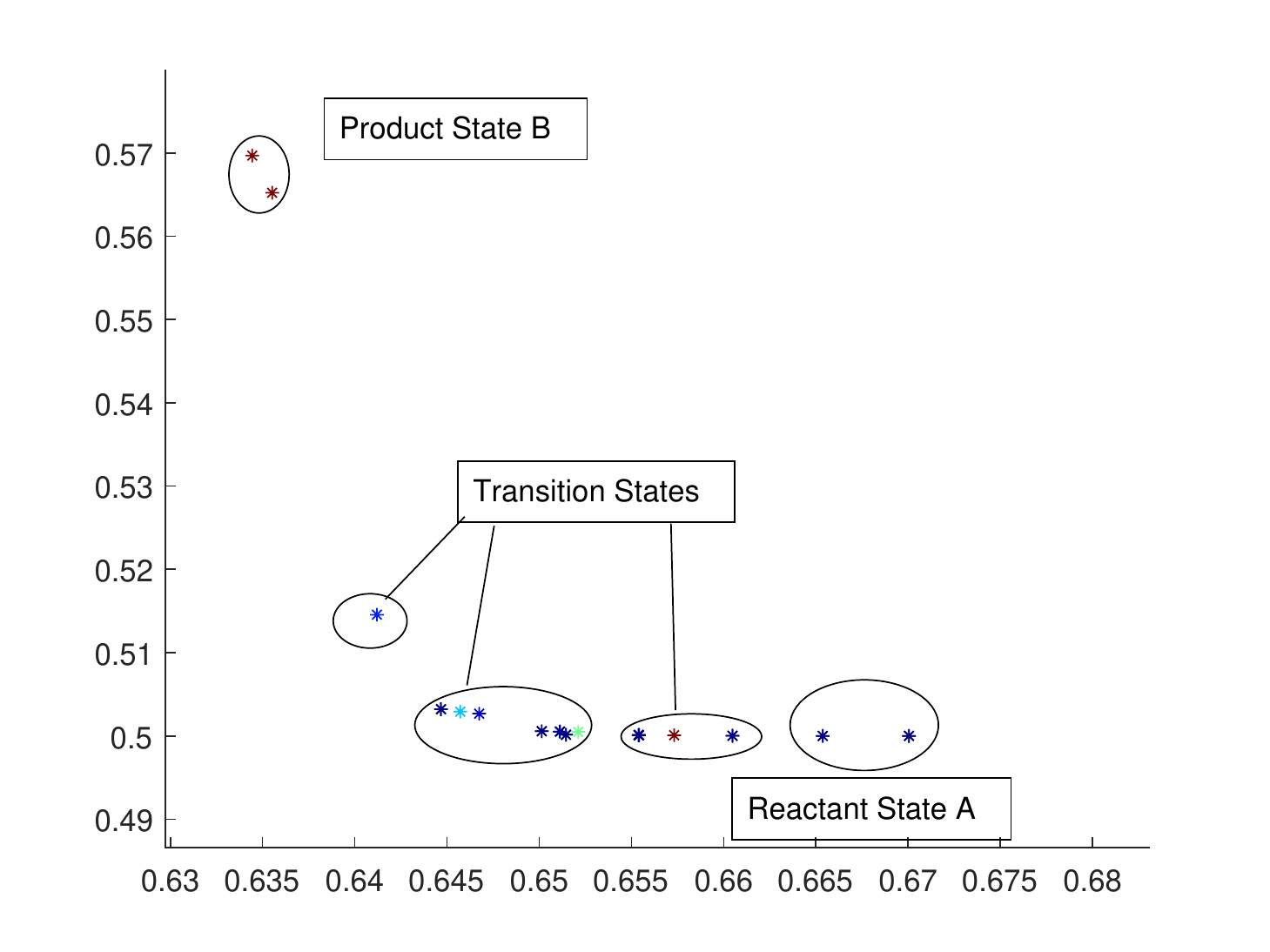}%
\caption{Feature representations of nodes with high similarity to the initial state $A=(0,0,0)$ for the virus 
propagation example. \label{fig:virusgroups}}%
\end{center}
\end{figure}

\begin{figure}[h] 
\begin{center}
\includegraphics[scale = 0.6]{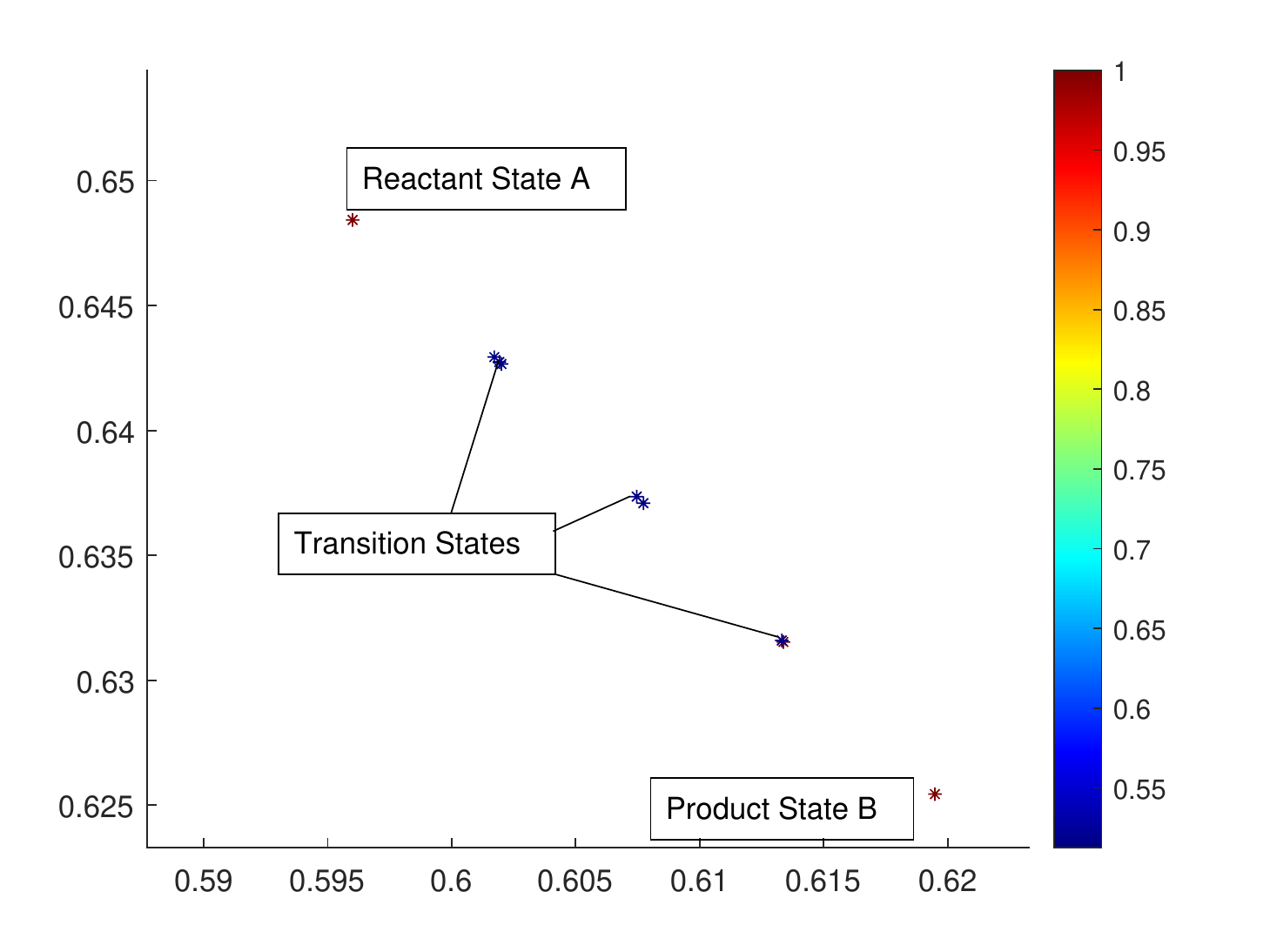}%
\caption{Feature representations for the virus propagation example, with a logarithmic scale on the $x$-axis.
\label{fig:logvirusgroups}}%
\end{center}
\end{figure}

The cellular concentrations of these proteins are controlled by six reactions: producing \emph{tem} from
\emph{gen}, using \emph{tem} as a catalyst to produce \emph{struct} and \emph{gen}, degradation of 
\emph{tem} and \emph{struct}, and propagation of the virus using \emph{struct} and \emph{gen} such that
\begin{align*}
&[gen] \xrightarrow{k_1} [tem]\ , \qquad
[tem] \xrightarrow{k_2} \emptyset,\\
&[tem] \xrightarrow{k_3} [gen]\ , \qquad
[gen]+[struct] \xrightarrow{k_4}\emptyset, \\
&[tem] \xrightarrow{k_5} [struct]\ , \quad
[struct] \xrightarrow{k_6} \emptyset,
\end{align*}
where we adopt $k_1=0.25, k_2=0.25, k_3=1.0, k_4=7.5 \times 10^{-6}, k_5=1000, \text{ and } k_6=1.99$. From 
\cite{stochvirus} we know there are two steady states for this system: an unstable trivial solution at 
$tem=gen=struct=0$,  and a stable steady state, which with the given parameter values is 
$tem=30, gen=100, struct=12000$. 
Ignoring the stochastic effects, as a macro scale limit, this system can be modeled by the ODEs:
\begin{align*}
\frac{d[tem]}{dt}&=k_1[gen]-k_2[tem],\\
\frac{d[gen]}{dt}&=k_3[tem]-k_4[gen][struct]-k_1[gen], \\
\frac{d[struct]}{dt}&=k_5[tem]-k_6[struct]-k_4[gen][struct]. 
\end{align*}
For a stochastic perspective, we model the system by a Markov process, where we follow the First Reactions 
method of the Gillespie algorithm\cite{gillespie}, which at each time step assumes that the reaction that could take 
place in the shortest amount of time will occur next. 



From the upper left to the lower right, these clusters contain the representations for nodes near approximately 
$(15 ,0 ,1000 )$, $( 30, 66.7, 4000)$, and $(30 , 100, 8000)$. Direct simulations confirm that these are transition 
states: reactive trajectories will tend to pass through nodes in each of these clusters. 
Applying the feature learning method to this system produces a two-dimensional feature representation for each
three-dimensional node $(tem, gen, struct)$. Here, we take state $A$ to be the origin and $B=(30, 100, 12000)$. 
Figure 7 plots the two-dimensional feature representations for this system. In this figure, the nodes with higher 
similarities to the initial state (represented in red) and those with lower such similarities (represented in blue) appear 
in separate clusters. Figure 8 plots these feature representations with a logarithmic scale on the $x$-axis to better
display the clusters. 
 
Five distinct clusters are discernible in Figure\ref{fig:logvirusgroups}: the red node in the top left is the feature 
representation for the reactant state, the red node at the bottom right is the representation of the product state, and 
the three clusters between each represent a transition state that reactive trajectories pass through while traveling 
from $A$ to $B$. From the upper left to the lower right, these three clusters contain the representations for nodes
near approximately $(15 ,0 ,1000 )$, $( 30, 66.7, 4000)$, and $(30 , 100, 8000)$. Direct simulations confirm that 
these are transition states: reactive trajectories will tend to pass through nodes in each of these clusters.
 
 \subsection{E. Coli Sigma-32 Heat Response Circuit}
Finally we will address a higher-dimension example, to illustrate the effectiveness of the method on reducing the 
dimension of a graph. Consider the E. coli $\sigma$-32 stress circuit, a network of regulatory pathways controlling
the $\sigma$-32 protein, which is essential in the E. coli response to heat shock \cite{sigma32}.
The systems consists of 10 processes:
\begin{align*}
&FtsH \xrightarrow{k_1}\emptyset, 
&&E\sigma^{32} \xrightarrow{k_2}  FtsH , \\
&GroEL \xrightarrow{k_3}\emptyset , 
&&E\sigma^{32}\xrightarrow{k_4}  GroEL,\\
&\sigma^{32} + J_{comp}  \xrightarrow{k_5}\sigma^{32}\text{-}J_{comp}, 
&&\sigma^{32}\text{-}J_{comp} \xrightarrow{k_6}  J_{comp} + \sigma^{32}, \\
&\sigma^{32}\text{-}J_{comp} \xrightarrow{k_7}  FtsH ,
&&\emptyset\xrightarrow{k_8}  \sigma^{32}, \\
&E+\sigma^{32} \xrightarrow{k_9} E\text{-}\sigma^{32}, 
&&E\text{-}\sigma^{32} \xrightarrow{k_10} E+\sigma^{32}.
\end{align*}
In this context, FtsH and GroEL are stress response proteins, E is a holoenzyme, and $J_{comp}$ (or J-complex)
represents several chaperone proteins that are lumped as a simplification. Then $E\text{-}\sigma^{32}$ denotes the
protein complex formed when E binds to $\sigma^{32}$, which catalyzes downstream synthesis reactions. 
$\sigma^{32}$ is a product of translation, which we assume to occur at a rate corresponding to $k_8=0.007$. 
It can also associate and dissociate of the J-complex. In accordance with \cite{sigma32} the other rate constants 
are taken to be: $k_1=7.4 \times 10^{-11} , k_2= 4.41 \times 10^{6}, k_3=1.80 \times 10^{-8} , k_4=5.69 \times 
10^6 , k_5=3.27 \times 10^{5} , k_6= 4.4 \times 10^{-4}, k_7=1.28 \times 10^{3} ,  k_9= 0.7, $ and $k_{10}= 0.13$.

For the above reaction rates, this system has a metastable states where the concentrations of FtsH, GroEL, 
$\sigma^{32}$, and $J_{comp}$ are approximately $600$, and the concentrations of the other three species are 
approximately $1500$. Another metastable state occurs when  the concentrations of FtsH, GroEL,$\sigma^{32}$, 
and $J_{comp}$ are approximately $800$. Here we take the former to be initial state $A$ and the latter to be the 
product state $B$. Following a procedure similar to that in the above virus propagation example, three-dimensional
representations are generated for each node. The plot of these node representations is given in Figure 9. The color
scale in this figure is determined by the value of the similarity function between each node and the node $A$, with 
the most similar nodes shown in red. Nodes with lower similarities have been omitted for clarity. Four distinct 
clusters can be seen in Figure 9.  

\begin{figure}[H] 
\begin{center}
\includegraphics[scale = 0.8]{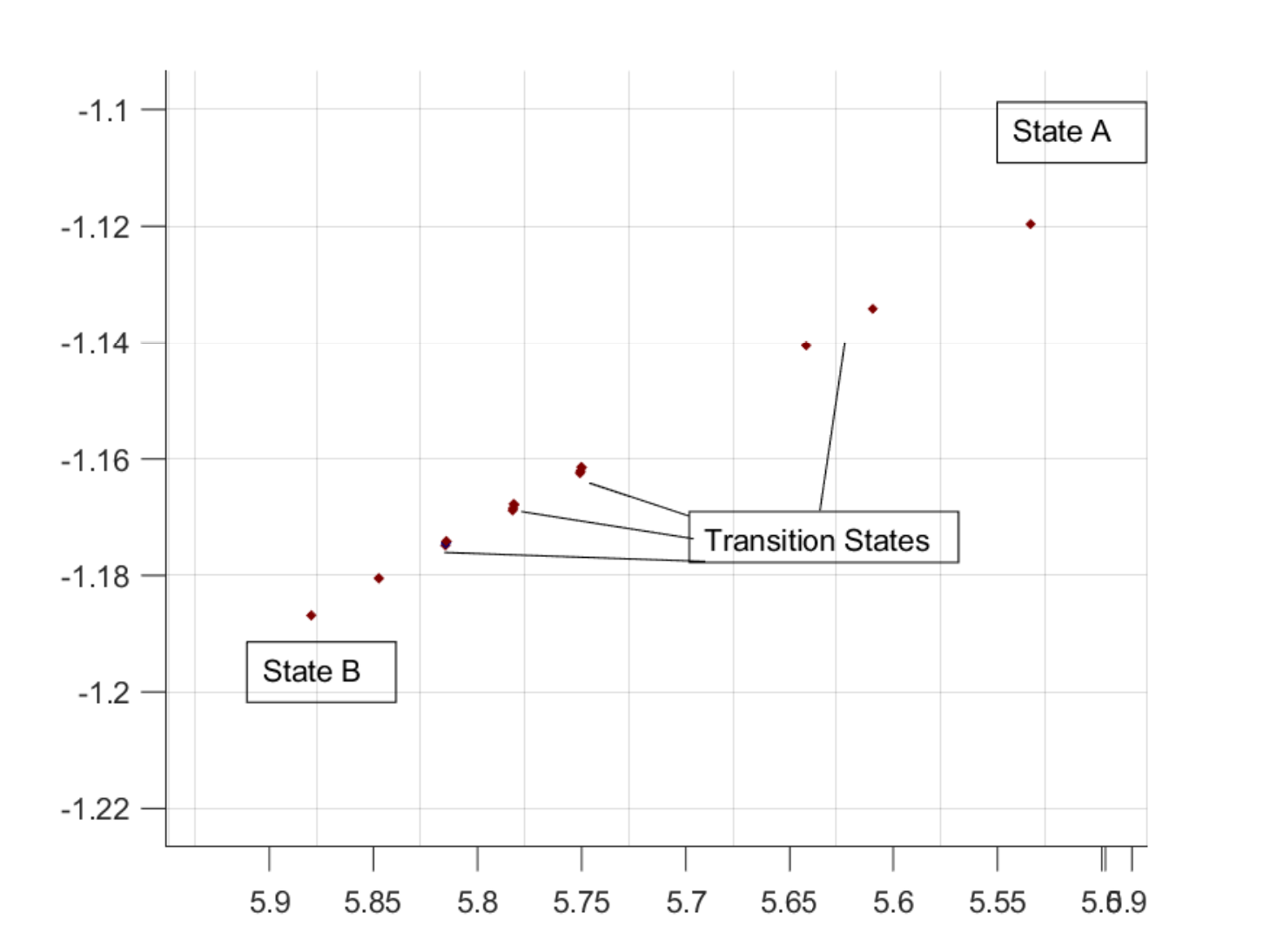}%
\caption{Feature representations of nodes with high similarity to the initial state, when initial state is taken to be $A=(600,600,600,600,1500,1500,1500)$, for the E. Coli heat response circuit model.\label{fig:sigma32}}%
\end{center}
\end{figure}

\section{Conclusions}
In this paper, we presented a method for analyzing metastable chemical reaction systems via feature learning on 
directed networks using random walk sampling. We have shown how this method may be used to identify transition
states of Markov jump processes by interpreting such processes in terms of directed networks and taking 
advantage of Transition Path Theory. We have illustrated the efficacy of this method through several 
low-dimensional examples involving various energetic and entropic barriers. As noted above, for more complex, 
realistic examples, the method can be used for dimensional reduction, enabling the extraction and analysis of 
high-dimensional graph information. Further work will be necessary to fully understand this approach from a 
probabilistic standpoint. Another direction for future work is improving efficiency and development of faster 
algorithms. However, this method still forms the groundwork for a potential new way of analyzing directed 
networks and jump processes, particularly in the context of chemical kinetics.


\begin{thebibliography}{10}

\bibitem{ahmedetal}
\newblock A.~Ahmed, N.~Shervashidze, S.~Narayanamurthy, V.~Josifovski and
  A.~Smola,
\newblock Distributed large-scale natural graph factorization,
\newblock in \emph{WWW 2013 - Proceedings of the 22nd International Conference
  on World Wide Web}, 2013,
\newblock 37--48.

\bibitem{bionetworks}
\newblock R.~Allen, P.~Warren and P.~Wolde,
\newblock Sampling rare switching events in biochemical networks,
\newblock \emph{Physical review letters}, \textbf{94} (2005), 018104.

\bibitem{lapeigs}
\newblock M.~Belkin and P.~Niyogi,
\newblock Laplacian eigenmaps and spectral techniques for embedding and
  clustering,
\newblock \emph{Advances in Neural Information Processing System}, \textbf{14}.

\bibitem{b:Wigner38}
E. Wigner, \emph{The transition state method}, Trans. of the Faraday Society, 34, 29-41, 1938.

\bibitem{b:Bolhuis02}
P. Bolhuis, D. Chandler, C. Dellago and P. Geissler,
\emph{Transition path sampling: Throwing ropes over rough mountain passes in the dark}, 
Ann. Rev. Phys. Chem., 53, 291-318, 2002.

\bibitem{b:EV06}
W. E and E. Vanden-Eijnden, 
\emph{Towards a theory of transition paths}, 
J. Stat. Phys., 123, 503-523, 2006.

\bibitem{survey2}
\newblock H.~Cai, V.~Zheng and K.~Chang,
\newblock A comprehensive survey of graph embedding: Problems, techniques and
  applications,
\newblock \emph{IEEE Transactions on Knowledge and Data Engineering}.

\bibitem{grarep}
\newblock S.~Cao, W.~Lu and Q.~Xu,
\newblock Grarep,
\newblock in \emph{CIKM '15: Proceedings of the 24th ACM International on
  Conference on Information and Knowledge Management}, 2015,
\newblock 891--900.

\bibitem{directedreps}
\newblock M.~Chen, Q.~Yang and X.~Tang,
\newblock Directed graph embedding,
\newblock in \emph{IJCAI International Joint Conference on Artificial
  Intelligence}, 2007,
\newblock 2707--2712.

\bibitem{cheegerineq}
\newblock F.~Chung,
\newblock Laplacians and the cheeger inequality for directed graphs,
\newblock \emph{Annals of Combinatorics}, \textbf{9} (2005), 1--19.

\bibitem{du}
\newblock J.~Du and D.~Liu,
\newblock Transition states of stochastic chemical reaction networks,
\newblock \emph{Comm. Comp. Phys.}, to appear.

\bibitem{handbook}
\newblock C.~Gardiner, U.~Bhat, D.~Stoyan, D.~Daley, Y.~Kutoyants and B.~Rao,
\newblock Handbook of stochastic methods for physics, chemistry and the natural
  sciences.,
\newblock \emph{Biometrics}, \textbf{42} (1986), 226.

\bibitem{gillespie}
\newblock D.~Gillespie,
\newblock Approximate accelerated stochastic simulation of chemically reacting
  systems,
\newblock \emph{Journal of Chemical Physics}, \textbf{115} (2001), 1716--1733.

\bibitem{survey1}
\newblock P.~Goyal and E.~Ferrara,
\newblock Graph embedding techniques, applications, and performance: A survey,
\newblock \emph{Knowledge-Based Systems}.

\bibitem{node2vec}
\newblock A.~Grover and J.~Leskovec,
\newblock node2vec: Scalable feature learning for networks,
\newblock in \emph{KDD : proceedings. International Conference on Knowledge
  Discovery \& Data Mining}, vol. 2016, 2016,
\newblock 855--864.

\bibitem{tptexamples}
\newblock P.~Metzner, C.~Schutte and E.~Vanden-Eijnden,
\newblock Illustration of transition path theory on a collection of simple
  examples,
\newblock \emph{The Journal of chemical physics}, \textbf{125} (2006), 084110.

\bibitem{tpt}
\newblock P.~Metzner, C.~Schutte and E.~Vanden-Eijnden,
\newblock Transition path theory for markov jump processes,
\newblock \emph{Mult. Mod. Sim.}, \textbf{7} (2008), 1192--1219.

\bibitem{negsampling}
\newblock T.~Mikolov, I.~Sutskever, K.~Chen, G.~Corrado and J.~Dean,
\newblock Distributed representations of words and phrases and their
  compositionality,
\newblock \emph{Advances in Neural Information Processing Systems},
  \textbf{26}.

\bibitem{deepwalk}
\newblock B.~Perozzi, R.~Al-Rfou and S.~Skiena,
\newblock Deepwalk: Online learning of social representations,
\newblock \emph{Proceedings of the ACM SIGKDD International Conference on
  Knowledge Discovery and Data Mining}.

\bibitem{lle}
\newblock S.~Roweis and L.~Saul,
\newblock Nonlinear dimensionality reduction by locally linear embedding,
\newblock \emph{Science (New York, N.Y.)}, \textbf{290} (2001), 2323--6.

\bibitem{gnn}
\newblock F.~{Scarselli}, M.~{Gori}, A.~C. {Tsoi}, M.~{Hagenbuchner} and
  G.~{Monfardini},
\newblock The graph neural network model,
\newblock \emph{IEEE Transactions on Neural Networks}, \textbf{20} (2009),
  61--80.

\bibitem{sigma32}
\newblock R.~Srivastava, M.~S. Peterson and W.~E. Bentley,
\newblock Stochastic kinetic analysis of the escherichia coli stress circuit
  using $\sigma$32-targeted antisense,
\newblock \emph{Biotechnology and Bioengineering}, \textbf{75} (2001),
  120--129,
\newblock
  \urlprefix\url{https://onlinelibrary.wiley.com/doi/abs/10.1002/bit.1171}.

\bibitem{stochvirus}
\newblock R.~Srivastava, L.~You, J.~Summers and J.~Yin,
\newblock Stochastic vs. deterministic modeling of intracellular viral
  kinetics,
\newblock \emph{Journal of theoretical biology}, \textbf{218} (2002), 309--21.

\bibitem{zhou}
\newblock D.~Zhou, J.~Huang and B.~Schoelkopf,
\newblock Learning from labeled and unlabeled data on a directed graph,
\newblock \emph{Framework}.

\end{thebibliography}
\providecommand{\href}[2]{#2}
\providecommand{\arxiv}[1]{\href{http://arxiv.org/abs/#1}{arXiv:#1}}
\providecommand{\url}[1]{\texttt{#1}}
\providecommand{\urlprefix}{URL }

\end{document}